\documentclass{amsart}
\usepackage{amsfonts}
\usepackage{amssymb}
\usepackage{color}

\newtheorem{thm}{Theorem}

\newtheorem{cor}[thm]{Corollary}
\newtheorem{lem}[thm]{Lemma}
\newtheorem{prop}[thm]{Propositon}
\newtheorem{rem}[thm]{Remark}
\newtheorem{exm}[thm]{{\it Example}}

 \newcommand{\lin}{\textrm{lin}}
\newcommand{\Real}{\mathbb{R}}
\newcommand{\Comp}{\mathbb{C}}
\newcommand{\Nat}{\mathbb{N}}
\newcommand{\Zzz}{{\mathbb{Z}}}
\newcommand{\eps}{\varepsilon}
\newcommand{\dts}{,\dots,}

\newcommand{\sbs}{\subseteq}
\newcommand{\h}{\mathcal{H}}

\newcommand{\Dd}{\mathcal{D}}

\newcommand{\Ssim}{S/_\sim}

\newcommand{\To}{ \Rightarrow}
\newcommand{\BOP}{\bold{B}}

\newcommand{\Pp}[1]{\frak{P}(#1)}

\newcommand{\Aa}{\mathcal{A}}
\newcommand{\Ll}{\mathcal{L}}
\newcommand{\M}[1]{\frak{M}(#1)}

\newcommand{\cl}{\overline}

\newcommand{\rest}[1]{\!\!\mid_{#1}}

\newcommand{\set}[1]{\left\{#1\right\}}
\newcommand{\seq}[1]{\left<#1\right>}
\newcommand{\norm}[1]{\left\Vert#1\right\Vert}

\begin{document}

   \title[Pontryagin space structure in RKHS's]{Pontryagin space structure in
   reproducing kernel Hilbert spaces over $*$--semigroups}
   \author[F.H. Szafraniec, M. Wojtylak]{ Franciszek Hugon Szafraniec \and Micha\l{} Wojtylak }

   \address{F. H. Szafraniec,  Instytut        Matematyki,         Uniwersytet
   Jagiello\'nski, ul. \L ojasiewicza 6, 30 348 Krak\'ow, Poland
    M. Wojtylak, Instytut        Matematyki,         Uniwersytet
   Jagiello\'nski, ul. \L ojasiewicza 6, 30 348 Krak\'ow, Poland\\ VU University Amsterdam,
 Department of Mathematics,
             Faculty of Exact Sciences,
             De Boelelaan 1081 a, 1081 HV Amsterdam\\
            }

            \email{umszafra@cyf-kr.edu.pl}
\email{michal.wojtylak@gmail.com}

   \thanks{The first author was  supported  by the MNiSzW grant
N201 026 32/1350. He also would like to acknowledge an assistance of
the EU Sixth Framework Programme for the Transfer of Knowledge
``Operator theory methods for differential equations'' (TODEQ) \#
MTKD-CT-2005-030042.}


   \begin{abstract}
The geometry of spaces with indefinite inner product, known also as
Krein spaces, is  a basic tool for developing Operator Theory
therein. In the present paper we establish a link between this
geometry  and the algebraic theory of $*$-semigroups.  It goes via
the positive definite functions and related to them reproducing
kernel Hilbert spaces. Our concern is in describing properties of
elements of the semigroup which determine shift operators which
serve as Pontryagin fundamental symmetries.
   \end{abstract}
   \maketitle

\section*{Introduction}


There are two ways of looking at $*$--semigroups and positive
definite functions defined on them. The first consists in
 intense analysis of the algebraic structure of a semigroup so as
to establish conditions on it, which ensure prescribed properties to hold for \underbar{any}
positive definite functions. One of the
 properties frequently considered is representing \underbar{each}
of the positive definite functions as moments of a positive
measure. This attitude has been successfully undertaken by
Bisgaard resulting in a considerable number of papers, see in
particular
\cite{biss-scmat,biss-two,biss-sep,biss-ext,biss-CC,biss-fact} and
references therein, some of them we are going to exploit here.
 The other way is to determine \underbar{which} of the positive
definite functions posses desired properties; a typical example
within this category is to detect multidimensional moment
sequences.

  In the present paper we are going to develop the first thread
putting forward the following problem. Impose necessary and sufficient conditions on a
distinguished element $u$ of a $*$--semigroup $S$ to generate a Pontryagin fundamental
symmetry of \underbar{any} reproducing kernel Hilbert space over the semigroup in question;
the precise formulation is exposed as ({\tt P}), p. \pageref{a}. Surprisingly, our problem
has found a simple algebraic solution. In the context of $*$-separative commutative
semigroups the condition on $u$ is: $u=u^*$, $u+u=0$ and $u+s\neq s$ for only a finite
number of $s\in S$ (see Proposition \ref{krein}, Theorem \ref{main}).

Let us mention that $*$-semigroups and positive definite functions on them have been
originated by Sz.-Nagy in his famous Appendix \cite{nagy}. He uses the reproducing kernel
Hilbert space factorization to prove his "general dilation theorem" for operator valued
functions. Since then the RKHS technique has been used from time to time for proving results
with dilation flavour behind the screen. We are going to follow suit here.


\section{Shift operators connected with positive definite
functions -- formulation of the problem}

By a $*$-semigroup we understand a commutative semigroup with an involution, not necessarily
having the neutral element. The involution is always denoted by the symbol ``$*$'' and the
semigroup operation is always written in an additive way. In the case when the
$*$--semigroup $S$ has a neutral element $0$ we say that $\phi:S\to\Comp$ is {\it positive
definite} (we write $\phi\in\Pp S$) if for every $N\in\Nat:=\set{0,1,\dots}$ and every
$s_0\dts s_N\in S$, $\xi_0\dts \xi_N\in\Comp$ we have $\sum_{i,j=0}^N\xi_i\bar\xi_j
\phi(s_j^*+s_i)\geq 0$. With such $\phi$ we link a reproducing kernel Hilbert space
$\h^\phi\sbs \Comp^S$ with a reproducing kernel defined by $K^\phi(s,t):=\phi(t^*+s)$,
$t,s\in S$ (see e.g. ~\cite[p.81]{BChR}). For $s\in S$ we set $K^\phi_s:=K(\cdot,s)$ ($s\in
S$), it is known that $K^\phi_s\in\h^\phi$ and $f(s)=\seq{f,K^\phi_s}$ for every
$f\in\h^\phi$. The set $\Dd^\phi:=\lin\set{K^\phi_s:s\in S}$ is dense in $\h^\phi$.

For an element  $u\in S$ we define {\it the  shift operator}
 (sometimes called also the translation operator) $
A(u,\phi):\Dd^\phi\to\Dd^\phi$, by
 $$
 A(u,\phi)\left(\sum_j\xi_jK_{s_j}^\phi\right):=\sum_j\xi_jK_{s_j+u}^\phi.
 $$
 It can be shown that $A( u,\phi)$ is well defined, linear and
closable (see e.g. \cite[Proposition, p.253]{szaf-bdd},
\cite[p.90]{BChR}). Our main object will be the closure of the
above operator, denoted below by $u_\phi$. It is a matter of
simple verification that $\Dd^\phi$ is contained in the domain of
$(u_\phi)^*$ (the adjoint of the operator $u_\phi$) and
$(u_\phi)^*f=(u^*)_\phi f$ for $f\in\Dd^\phi$.



Under the above circumstances we can state our problem as follows.

    \begin{enumerate} \label{a}
   \item[({\tt P})] Provide necessary and sufficient conditions on
   the element $u\in S$ for the operator $u_\phi$ to be a
   fundamental symmetry of a Pontryagin space, i.e. to satisfy
$$
u_\phi=u_\phi^*,\quad u_\phi^2=I_{\mathcal{H}^\phi},\quad
\dim\ker(u_\phi+I_{\mathcal{H}^\phi})<\infty
$$
 \underbar{for every} $\phi\in\Pp S$.

   \end{enumerate}
Obviously, the
 condition $u_\phi=u_\phi^*$ together with $u_\phi^2=I_{\mathcal{H}^\phi}$ imply that $u_\phi$ must be bounded
on $\h$.

We continue with some basic definitions and notations concerning $*$-semigroups. Let $S$ and
$T$ be  $*$-semigroups, a mapping $\chi:S\to T$ is called $*$-{\it homomorphism} if
$\chi(s+t)=\chi(s)+\chi(t)$ for all $s,t\in S$, $\chi(s^*)=(\chi(s))^*$ for all $s\in S$. A
{\it character} on $S$ is a nonzero $*$-homomorphism $\chi:S\to\Comp$ where the latter set
is understood as a semigroup with multiplication as the operation and the conjugation as
involution. It is obvious that if $S$ has a neutral element $0$ then $\chi(0)=1$ for all
$\chi\in S^*$. Let $\Aa(S^*)$ be the least $\sigma$-algebra of subsets of $S^*$ rendering
measurable all the functions
\begin{equation}\label{shat}
\hat s:S^*\ni\chi\mapsto\chi(s)\in\Comp,\quad s\in S.
 \end{equation}
  If $\mu$ is a
positive measure on $S^*$ such that all the functions listed in
(\ref{shat}) are square-integrable we define a function
$\Ll(\mu):S\to\Comp$ by
$\Ll(\mu)(s):=\int_{S^*}\sigma(s)d\mu(\sigma)$ ($s\in S$). We call
$\phi:S\to\Comp$ a {\it moment function} ($\phi\in\M S $) if
$\phi=\Ll(\mu)$ for some measure $\mu$ on $S^*$. It is easy to
verify that $\M S \sbs\Pp S$, if the latter inclusion is an equality
then we call $S$ {\it semiperfect}. For examples of semiperfect and
non-semiperfect $*$-semigroups and more general concepts of
semiperfectness see ~\cite{biss-scmat}.

We call $S$ {\it $*$-separative} if the characters separate points in $S$. The {\it greatest
$*$-homomorphic $*$-separative image of $S$} is the semigroup $\Ssim$, where the equivalence
relation $\sim$ on $S$ is defined by the condition that $s\sim t$ if and only if
$\sigma(s)=\sigma(t)$ for all $\sigma\in S^*$; addition and involution in $\Ssim$ are those
that make the quotient mapping a $*$-homomorphism. The elements of $\Ssim$ will be denoted
as equivalence classes $[s]$ ($s\in S$). If $S$ has a zero, then we will use the symbol $0$
for the neutral element of both $S$ and $\Ssim$, instead of using $[0]$. Since every
character on $S$ generates a character on $\Ssim$, the latter semigroup is in fact
$*$-separative. The following simple proposition gives answer to the question when the
operator $ u_\phi$ defined above is a fundamental symmetry of a Krein space, i.e. when $(
u_\phi)^2= u_\phi$ and $ u_\phi=( u_\phi)^*$.

\begin{prop}\label{krein}
Let $S$ be a commutative $*$-semigroup with zero. For each element $u\in S$ the following
conditions are equivalent:
\begin{itemize}
\item[(i)]{$[2u]= 0$ and $[u]=[u^*]$;}
\item[(ii)]{for every $\phi\in\M S $ we have $(u_\phi)^2=
I_{\h^\phi}$ and $(u_\phi)^*= u_\phi$;} \item[(iii)]{for every
$\phi\in\Pp S$ we have $(u_\phi)^2 = I_{\h^\phi}$ and $(u_\phi)^*=
u_\phi$;}
\end{itemize}
Moreover, {\rm (i)} implies that  $ u_\phi$ is a bounded,  selfadjoint operator on $\h^\phi$ for every
$\phi\in\Pp S$.
\end{prop}

\begin{proof}
(i)$\To$(iii) Suppose  that (i) holds and let $\phi\in\Pp S$. For every $\sigma\in S^*$ we have
$\sigma(u)=\sigma(u^*)$ and $\sigma(2u)=\sigma(0)=1$. Consequently, for every $\sigma\in S^*$ and every $s\in
S$
\begin{equation}\label{sigma2}
\sigma(s+u)=\sigma(s+u^*),\quad\sigma(t+2u)=\sigma(t),\quad \sigma(s^*+u^*+u+s)=\sigma(s^*+s).
\end{equation}
 By \cite[Thm.2]{biss-sep}  we have
\begin{equation}\label{sigma3-1}
\phi(s+u)=\phi(s+u^*),\quad s\in S,
\end{equation}
\begin{equation}\label{sigma3-2}
\phi(s+2u)=\phi(s),\quad s\in S,
\end{equation}
\begin{equation}\label{sigma3-3}
 \phi(s^*+u^*+u+s)=\phi(s^*+s),\quad s\in S,
\end{equation}
The last of these three equalities, together with \cite[Cor.1]{szaf-bdd}, implies that the
operator $u_\phi$ is in $\BOP(\h^\phi)$. It is also selfadjoint, since for $f=\sum_i \xi_i
K_{s_i}^\phi\in\Dd^\phi$ we have
$$
 \seq{u_\phi f,f}=
 \sum_{i,j}\xi_i\bar\xi_j\seq{K_{s_i+u}^\phi,K_{s_j}^\phi}=\sum_{i,j}\xi_i\bar\xi_j\phi(s_i+u+s_j^*)
 \stackrel{(\ref{sigma3-1})}=
 $$ $$
 \sum_{i,j}\xi_i\bar\xi_j\phi(s_i+u^*+s_j^*)=
 \seq{f, u_\phi f}.
 $$
 The fact  that $(u_\phi)^2=I_\h^\phi$ can be obtained
similarly as selfadjointness of $u_\phi$, with the use of
(\ref{sigma3-2}) instead of (\ref{sigma3-1}). This finishes the
proof of (iii) and of the `Moreover' part of the proposition.

The implication (iii)$\To$(ii) is trivial. The proof (ii)$\To$(i) goes by contraposition.
Suppose first that $[2u]\neq 0$, i.e. there exists a character $\sigma$ such that
$\sigma(u)^2=\sigma(2u)\neq\sigma(0)=1$. We put $\phi:=\Ll(\delta_\sigma)$, where
$\delta_\sigma$ stands for the Dirack measure on $S^*$ concentrated in $\sigma$. Let
$\seq{\,\cdot\,,-}$ denote the scalar product on $\h^\phi$. Observe that
$$
\seq{K_0^\phi,K_0^\phi}=\phi(0)=\sigma(0)\neq\sigma(2u)=\phi(2u)=\seq{K_{2u}^\phi,K_0^\phi}=\seq{(
u_\phi)^2K_0^\phi,K_0^\phi}.
$$
In consequence, $I_{\h^\phi}\neq ( u_\phi)^2$. Similarly, if $\tau(u)\neq\tau(u^*)$ for some
$\tau\in S^*$ then for $\psi:=\Ll(\delta_\tau)$ the operator $u_\psi$ is not symmetric in
$\h^\psi$.

\end{proof}

Let us consider now a situation when $S$ and $T$ are  $*$-semigroups
with zeros and $h$ is a $*$-homomorphism from $S$ into $T$
satisfying $h(0)=0$. Note that if an element $u\in S$ is such that
$[2u]=0$, $[u]=[u^*]$ then $[2h(u)]=0$, $[h(u)]=[h(u)^*]$. This
comes from the fact that for every character $\sigma$ on $T$ the
function $\sigma\circ h$ is a character on $S$. Observe also that
$\phi\circ h\in\Pp S$ for every $\phi\in\Pp T$.

\begin{prop}\label{ST}
Assume that $S$,  $T$ and $h$ are as above and that $h$ is
additionally onto. Let $u\in S$ be such that $[2u]=0$, $[u]=[u^*]$
and let $\phi\in\Pp T$. Then the operators $u_{\phi\circ h}$ in
$\h^{\phi\circ h}$ and $h(u)_{\phi}$ in $\h^\phi$ are unitarily
equivalent.
\end{prop}

\begin{proof}
 Let $\seq{\cdot,-}_\phi$ and
$\seq{\cdot,-}_{\phi\circ h}$ denote the scalar products on $\h^\phi$ and $\h^{\phi\circ h}$
respectively. Since
   \begin{eqnarray*}
   \seq{\sum_{j=1}^N \xi_j K^{\phi\circ h} _{s_j} ,\sum_{j=1}^N
\xi_j K^{\phi\circ h} _{s_j}}_{\!\!\phi\circ h}
 &=& \sum_{i,j=1}^N \xi_i\bar\xi_j (\phi\circ h)(s_i+s_j^*) \\ =
\sum_{i,j=1}^N \xi_i\bar\xi_j \phi(h(s_i)+h(s_j^*))&=&
 \seq{\sum_{j=1}^N \xi_j K^{\phi} _{h(s_j)},\sum_{j=1}^N \xi_j
K^{\phi} _{h(s_j)}}_{\!\!\phi},
   \end{eqnarray*}
the condition   $ V( K^{\phi\circ h}_{s}):= K^{\phi}_{h(s)} $ ($s\in
S$) properly defines an isometry between $\h^{\phi\circ h}$ and
$\h^{\phi}$. Since $h$ is onto, the range of $V$ is dense in
$\h^\phi$, and so $V$ is a unitary operator. To finish the proof we
need to show that
\begin{equation}\label{VuuV}
h(u)_\phi Vf=V u_{\phi\circ h}f, \quad f\in\h^{\phi\circ h}.
\end{equation}
This can be easily verified for $f\in\Dd^{\phi\circ h}$. Since all
the operators appearing in (\ref{VuuV}) are bounded, the proof is
finished.
%
%

\end{proof}

Applying the above to the quotient semigroup $T=\Ssim$ and $h$ as
the quotient mapping, together with the fact from \cite{biss-fact}
that every $\phi\in\Pp S$ ($\phi\in\M S$) is of the form
$\phi=\psi\circ h$ for some $\psi\in\Pp{\Ssim}$ ($\psi\in\M\Ssim$)
gives  the following.

\begin{cor}\label{Nowyjeszcze}
Assume that $S$ is a $*$-semigroup with zero and $u\in S$ is such
that $[2u]=0$, $[u]=[u^*]$. Then
 \begin{eqnarray*}
 \dim\ker(u_\phi+I)&<&+\infty \textrm{ for every }\phi\in\Pp S\,\, (\phi\in\M S)\\
   \iff \dim\ker([u]_{\psi}+I)&<&+\infty\textrm{ for every
 }\psi\in\Pp\Ssim\,\, (\psi\in\M\Ssim)
 \end{eqnarray*}
\end{cor}
\section{Examples}

%
%
%

\begin{exm}\label{exm1}
Consider the semigroup $S=\Zzz_2\times\Nat$ with standard addition
and the identical involution. As usually (cf. \cite{BChR}) we
identify $S^*$  with $\set{-1,1}\times\Real$, note that  $S$ is
$*$-separative. The only nonzero element satisfying $u=u^*$, $2u=0$
is $u=(1,0)$. Let $\phi$ be a positive definite mapping, we will now
compute the eigenspaces of $ u_\phi$.    Since $S$ is semiperfect
(\cite{biss-two}), there exists a Borel measure $\mu$ on $S^*$ such
that $\phi=\Ll(\mu)$. Having our interpretation of characters in
mind we get
$$
\phi(x,n)=\int_\Real (-1)^xt^nd\mu_-(t)+\int_\Real t^nd\mu_+ (t),\quad
x\in\Zzz_2,\,n\in\Nat,
$$
with $\mu_\pm:=\mu\rest{{\set{\pm1}}\times\Real}$. Let us define the functions
$f_{x,n}:S^*\to S^*$ ($x\in\Zzz_2$, $n\in\Nat$) by
 $$
 f_{x,n}(\eps,t):=\eps^xt^x ,\quad \eps\in\set{-1,1},\,t\in\Real,\,x\in\Zzz_2,\,n\in\Nat,
 $$
and note that they all are square integrable. By $\mathcal P^\mu$ we define the closure in
$L^2(\mu)$ of the linear span of the functions $f_{x,n}$ ($x\in\Zzz_2$, $n\in\Nat$). The
formula
  \begin{equation}\label{Vdef}
 V(K^\phi_{x,n}):=f_{x,n},\qquad x\in\Zzz_2,\,n\in\Nat,
 \end{equation}
constitutes a unitary isomorphism between $\h^\phi$ and $\mathcal P^\mu$. The shift operator
$u_\phi$ is unitary equivalent (via $V$) to the following operator $M$
    $$
    (M f)(\eps,t):=\eps f (\eps,t),\qquad (\eps,t)\in\set{-1,1}\times\Real,\, f\in\mathcal{P}^\mu.
    $$
     It is not hard to see that
    $$
    \ker(M\pm I)=\set{f\in \mathcal{P}^\mu: f(\pm1,\cdot)=0\,\,\, \,
    \mu_\pm\textrm{--a.e}}.
    $$
Note that $\dim\ker(u_\phi\pm I)=\dim\ker(M\pm I)$  and the latter is finite dimensional if
and only if the support of $\mu_{\mp}$ is a finite set. In particular there exists a mapping
$\phi\in\Pp S=\M S $ such that $\dim\ker( u_\phi+I)=\infty$.
\end{exm}

At this point we present a useful for our purposes construction. Let $S$ and $T$ be two
disjoint $*$-semigroups (which only a formal restriction) and let $h:S\mapsto T$ be a
$*$-homomorphism.  We endow the set $S\cup T$ with the $*$-semigroup structure in the
following way. The addition on $S\cup T$, denoted by the same symbol `$+$', is defined by
 $$
 s+t:=\left\{\begin{array}{rcl}
    s+t & : & s,t\in S\textrm{ or } s,t\in T\\
    s+h(t) & : & s\in S,\, t\in T\\
    t+h(s) & : & t\in S,\, s\in T\\
        \end{array}\right.
    $$
The involution on $S\cup T$ (still denoted by `$*$') is such that its restriction to both
$S$ and $T$ is the original involution on $S$ and $T$, respectively.  We denote the above
constructed semigroup by $\mathrm{U}(S,T,h)$. The reader can easily check a general fact,
that if $S$ and $T$ are $*$-separative then $\mathrm{U}(S,T,h)$ is $*$-separative as well.

\begin{exm}\label{exm2} Consider a semigroup $S=\mathrm{U}(\Zzz_2,\Nat,h_0)$
where  $h_0(x)=0$ ($x\in \Zzz_2$).
The element $0_{\Zzz_2}$ is the neutral element of $S$. Take  $u:=1_{\Zzz_2}$, clearly
$u=u^*$ and $2u=0_{\Zzz_2}$. Let $\phi$ be any positive definite function on $S$ and suppose
that $f\in\ker( u_\phi+I)$. This means that for $n\in\Nat$
 $$
 f(n)=f(n+u)=\seq{f,K^\phi_{n+u}}=\seq{f,u_\phi K^\phi_n}=\seq{u_\phi f, K^\phi_n}=-f(n),
 $$
 hence $f\rest\Nat=0$. A similar calculation shows that $f(u)=-f(0_{\Zzz_2})$. Hence,
the eigenspace $\ker( e_\phi+I)$ is spanned by the single function
 $$
 f(s)=\left\{\begin{array}{rcl}
         0 &:& s\in\Nat\\
         1 & : & s=0_{\Zzz_2}\\
         -1 & : & s=1_{\Zzz_2}
         \end{array} \right.
         $$
 if $f\in\h^\phi$ or is trivial otherwise. Resuming, $\dim\ker(u_\phi+I)\leq 1$ for all
 positive definite $\phi$.
\end{exm}

\section{Main result}

\begin{thm}\label{main}
Let $S$ be a commutative  $*$-semigroup with zero and let $u\in S$ be such that $[2u]=0$ and
$[u]=[u^*]$. Then the following conditions are equivalent:
\begin{itemize}
\item[(i)] the set  $\set{[s]\in S/_\sim\colon [u+s]\neq[s]}$  is
finite\,\footnote{\ Cf. Remark \ref{stupidstyle}.};
 \item[(ii)]{  $\sigma(u)=-1$ for only finitely many $\sigma\in S^*$;}
 \item[(iii)]{$\dim\ker(u_\phi+I_{\h^\phi})<\infty$ for all $\phi\in\M S $;}
 \item[(iv)]{$\dim\ker(u_\phi+I_{\h^\phi})<\infty$ for all $\phi\in\Pp S$.}
\end{itemize}
Moreover, if {\rm (i)} holds then $\dim\ker(u_\phi+I_{\h^\phi})$ is
less or equal to the half of the number elements of the set
mentioned in {\rm (i)}.
\end{thm}

Let us stress here the solution of the main problem of the paper,
which follows from Proposition~\ref{krein} and Theorem~\ref{main} :
{\it the operator $ u_\phi$ is a fundamental symmetry of a
Pontryagin space for every $\phi\in\Pp S$ $(\!\!$ equivalently:
$\phi\in\M S$$)$  if and only if  $[2u]=0$, $[u]=[u^*]$ and the set
$\set{[s]\in\Ssim:[u+s]\neq[s]}$ is finite.}

\begin{rem}\label{stupidstyle}
Condition (ii) is equivalent to
\begin{itemize}
\item[(ii')]{\it $\sigma([u])=-1$ for only finitely many characters $\sigma$ on $\Ssim$,}
\end{itemize}
since the characters on $S$ and $\Ssim$ are in one-to-one
correspondence. This, together with Corollary \ref{Nowyjeszcze} (and
the remarks above it) allows us to reduce the proof of Theorem
\ref{main} to the case when $S$ is $*$-separative. However, note
that the conditions (i)--(iv) are \underbar{not} equivalent to the
following:
\begin{itemize}
\item[]{\it the set $\set{s\in S: u+s\neq s}$ is finite.}
\end{itemize}
Indeed, consider the following example. Let $S=\Zzz_4\times\Nat$
with the natural operation $+$ and the identical involution. The
greatest $*$-separative homomorphic image of $S$ is ($*$-isomorphic
with) $\Zzz_2\times\Nat$ and the quotient homomorphism maps
$u:=(2,0)$ to $(0,0)$. We have that $u+s\neq s$ for all $s\in S$ but
the set mentioned in (i) is empty. It remains on open problem if the
condition
\begin{itemize}
\item[]{\it the set $\set{s\in S: [u+s]\neq [s]}$ is finite.}
\end{itemize}
is equivalent to (i).
\end{rem}

Before the proof we introduce the notion of a $*$-archimedean
component of a semigroup. We call a $*$-semigroup $H$ (not
necessarily with 0) {\it $*$-archimedean} if for all $s,t\in H$
there exists $m\in\Nat\setminus\set0$ such that $m(s+s^*)\in t+ H$.
An {\it $*$-archimedean component} of a $*$-semigroup $S$ is a
maximal $*$-archimedean $*$-subsemigroup. Though $*$-archimedean
component is $*$-semigroup for itself it is possible for it not to
have the neutral element even if $S$ does. It can be shown that two
elements $s,t$ belong to the same $*$-archimedean component of $S$
if and only if $m(s+s^*)\in t+S$ and $n(t+t^*)\in s+S$ for some
$m,n\in\Nat\setminus\set0$. Furthermore, $S$ is the disjoint union
of its $*$-archimedean components, see \cite[Section 4.3]{clipres}
for the case of identical involution. The following Lemma was proven
in \cite{biss-ext} (Lemma 2), the proof for an arbitrary involution
requires minimal effort.
\begin{lem}\label{extchar}
If $H$ is a $*$-archimedean  component of a $*$-semigroup $S$ then every character on $H$ is
everywhere nonzero and extends to a character on $S$.
\end{lem}

 If $H$ and $K$ are two
$*$-archimedean components of $S$ then $H+K$ is contained in one single $*$-archimedean
component of $S$. If $(S_i)_{i\in I}$ is the family of all $*$-archimedean components of $S$
then we define the operation $+$ on $I$ by:
$$
i+j=k\textrm{ if and only if } S_i+S_j\sbs S_k,\quad i,j,k\in I.
$$
Since $S_i+S_i\sbs S_i$ for all $i\in I$ we have that $i+i=i$. Therefore $I$ is a
semilattice with the natural partial order given by the condition that $i\leq j$ if and only
if $i+j=j$. The following easy lemma is left as an exercise for the reader.

\begin{lem}\label{elat}
Let $S$ be a $*$-semigroup with zero and let $(S_i)_{i\in I}$ be the family of all
$*$-archimedean components  $S$. If $u\in S$ be such that $2u=0$, then $u$ belongs to the
same $*$-archimedean component as 0. In particular, $u+S_i\sbs S_i$ for all $i\in I$.
\end{lem}



\begin{proof}[Proof of Theorem ~\ref{main}]
As it was said in Remark \ref{stupidstyle} we may assume that $S$ is
$*$-separative. To prove (i)$\To$(iv) let us put
    $$
     U:=\set{ s\in S: u+ s= s}
    $$
 and suppose that the set $ S\setminus U$ contains only a
finite number $M$ of elements. Take $\phi\in\Pp S$. We show that
    $$
    \dim\ker( u_\phi+I)\leq M/2,
    $$
this will also prove the last statement of Theorem \ref{main}. Let us fix an arbitrary
$f\in\ker( u_\phi+I)$. We have
 $$
 f( s+ u)=\seq{f,u_\phi K_{s}^\phi}=\seq{u_\phi f,K_s^\phi}=-f( s),\quad s\in  S.
 $$
  This means
that $f\rest{ U}\equiv 0$. Observe that the relation
 $$
 aRb \Longleftrightarrow (a=u+b\textrm{ or }a=b)
 $$
is an equivalence relation on $ S\setminus U$ and that the each equivalence class contains
exactly two elements. Take any representees $s_1\dts s_{M/2}$ of the equivalence classes of
$R$. It is clear, that
 $$
 \ker( u_\phi+I)\sbs\lin\set{\delta_{s_i}-\delta_{s_i+u}:i=1\dts M/2}.
  $$
  Consequently $\dim\ker( u_\phi+I)\leq M/2$.\\

(iv)$\To$(iii) is obvious. (iii)$\To$(ii) Suppose that $\set{\sigma_n:n\in\Nat}$ is an
infinite set of characters satisfying $\sigma_n(u)=-1$ ($n\in\Nat$).  We define a measure
$\mu$ on $S^*$ by $\mu:=\sum_{n=0}^\infty 2^{-n} \delta_{\sigma_n}$ and we take a mapping
$\phi=\Ll(\mu)$. Note that for every $N\in\Nat$ and for every $s_0\dts s_N\in S$,
$\xi_0\dts\xi_N\in\Comp$ we have
    $$
    \bigg|\sum_{j=0}^N\xi_j\sigma_n(s_j)\bigg|^2\leq
    2^n\int_{S^*}\bigg|\sum_{j=0}^N\xi_j\sigma(s_j)\bigg|^2d\mu(\sigma)=
        2^n\sum_{i,j=0}^N\xi_i\bar\xi_j\int_{S^*}\sigma(s_i+s_j^*)d\mu(\sigma)=
     $$ $$
    =2^n\sum_{i,j=0}^N\xi_i\bar\xi_j\phi(s_i+s_j^*),\quad   n\in\Nat.
    $$
By the RKHS Test (\cite{test} and also \cite{test2}) we get that $\sigma_n\in\h^\phi$ for
$n=1,2,\dots$. Now observe that
    $$
    u_\phi(\sigma_n)(s)=\sigma_n(u+s)=-\sigma_n(s),\quad s\in S,n=1,2\dots.
    $$
Therefore $\sigma_n\in\ker( u_\phi+I)$, $n=1,2,\dots$. It remains to
show that the functions $\sigma_n$, $n\in\Nat$, are linearly
independent. But this results from the well known fact that all
characters are linearly independent\,\footnote{\ We can use the
following argument: For every $s\in S$ the function $\hat s$ is a
character on the dual semigroup $S^*$ and it is trivial that the
family $(\hat s)_{s\in S}$ separates elements of $S^*$. Proposition
2 of \cite{lacunae} (see also \cite[Proposition 6.1.8]{BChR}) says
that if $T$ is a semigroup and $C\sbs T^*$ separates points, then
the functions $\hat t\rest C$, $t\in T$ are linearly independent in
$\Comp^C$. We use this result for $T=S^*$ and $C=\set{\hat s:s\in
S}$, the functions $\hat\sigma\rest C$ can
be identified with characters on $S$.}.\\

(ii)$\To$(i) Suppose that the number of elements $ s\in  S$ satisfying $ u+ s\neq  s$ is
infinite. We show that there exists infinitely many characters $\sigma$ on $ S$ such that
$\sigma( u)=-1$.

Let $( S_i)_{i\in I}$ be the family of all $*$-archimedean components of $ S$. Set
 $$
 I_0:=\set{i\in I:  u+ s\neq  s\textrm{ for some } s\in S_i}.
 $$
Our assumption implies that\,\footnote{\ Remark \ref{eitheror} shows
that it is even equivalent to} either
\begin{equation}\label{Case1}
\textrm{ $ S_j$ is infinite for some $j\in I_0$}
\end{equation}
or
\begin{equation}\label{Case2}
I_0\textrm{ is infinite}.
\end{equation}

Let us first assume (\ref{Case1}). Take $s_0\in S_j$ such that
$u+s_0\neq s_0$. By Lemma \ref{elat} we have $ u+ s_0\in S_j$. The
$*$-semigroup $S_j$ is $*$-separative as a subsemigroup of a
$*$-separative semigroup $S$. Therefore, there exists a character
$\sigma_0$ on $ S_j$ such that $\sigma_0( u+ s_0)\neq\sigma_0( s_0).
$
 By
Lemma \ref{extchar} $\sigma_0$ extends to some character $\tilde{\sigma_0}$ on $ S$. Since $
u= u^*$  and $2 u=0$ we have  $\tilde{\sigma_0}( u)\in\set{-1,1}$. But
 $$
 \tilde{\sigma_0}( u)\sigma_0( s_0)=\sigma_0( u+ s_0)\neq\sigma_0( s_0).
    $$
Hence, $\tilde{\sigma_0}( u)=-1$, which means that $\sigma_0( u+ s_0)=-\sigma_0( s_0)$.

Denote by $A$ the set of all those characters $\sigma$ on $ S_j$ satisfying $\sigma( u+
s_0)=-\sigma( s_0)$. Since $\sigma_0$ is everywhere nonzero on $ S_j$ (Lemma \ref{extchar}),
the mapping
 $$
  S_j^*\ni\sigma\longmapsto \sigma_0\sigma\in S_j^*
 $$
is bijective. Moreover, it maps $A$ onto $ S_j^*\setminus A$.  Since $ S_j$ is infinite and
$*$-separative, $ S_j^*$ is infinite as well. Hence, $A$ is infinite. By Lemma
\ref{extchar}, there is an infinite number of characters $\sigma$ on $ S$ satisfying
$\sigma( u+ s_0)=-\sigma( s_0)$ and consequently $\sigma( u)=-1$.\\

Let us assume now (\ref{Case2}). For each $i\in I_0$ we take a
character $\sigma_i$ on $ S$ satisfying $\sigma_{i}( u)=-1$,
$\sigma_{i}( s)\neq 0$ for $s\in S_{i}$ (such a character exist by
repeating the proof from the previous case). We also define a family
of characters $\chi_{i}\in S^*$ (${i\in I_0}$)  by
    $$
    \chi_{i}( s)=\left\{\begin{array}{rl} 1 & \textrm{ if }  s\in S_j\textrm{ and }j\leq i\\
                                             0 & \textrm{
                                             otherwise}
                                             \end{array}\right.,\qquad  s\in  S,\,i\in I_0.
    $$
Finally, we put $\rho_i:=\sigma_i\chi_i$ ($i\in I_0$). By Lemma \ref{elat} we have that $ u$
is in the same $*$-archimedean component as 0, denote this component by $ S_{j_0}$. It is
clear that $j_0\leq i$ for all $i\in I$, therefore $\chi_i( u)=1$ and $\rho_i( u)=-1$ for
all $i\in I_0$. The only thing that lasts is to show that $\rho_i\neq\rho_j$ for $i\neq j$.
If $i\neq j$ then, by symmetry, we can assume that $j\nleq i$. Thus $\chi_i\rest{ S_j}=0$
and $\rho_i\rest{ S_j}=0$. But $\chi_j\rest{ S_j}\equiv1$ and $\sigma_j$ is, by definition,
everywhere nonzero on $ S_j$. Therefore
$\rho_i\rest{ S_j}\equiv0\neq\rho_j\rest{ S_j}$.\\

\end{proof}

\begin{rem}\label{eitheror}
The alternative of (\ref{Case1}) and (\ref{Case2}) in the proof of
(ii)$\To$(i) becomes more clear if we observe that $ u+ s\neq  s$
for \underbar{all} $ s\in S_i$, provided that $i\leq j\in I_0$.
Indeed, suppose that $ u+ s= s$ for some $ s\in S_i$, $ u+ s_0\neq
s_0$ for some $ s_0\in S_j$ and $i\leq j$. This gives us
 $$
 ( u+ s_0)+( s+ s_0)= u+ s+ s_0+ s_0= s_0+( s+ s_0).
 $$
  By Lemma \ref{elat} we have
$ u+ s_0\in S_j$. Moreover, $ s+ s_0\in S_j$ because $i\leq j$. The
semigroup $ S_j$ is cancellative as a $*$-archimedean component of a
$*$-separable group (see \cite[p.63]{biss-CC}). This gives us $ u+
s_0= s_0$, contradiction. The example below shows that both
(\ref{Case1}) and (\ref{Case2}) are possible.

\end{rem}

\begin{exm}
Let $S=\Zzz_2\times\Nat$, with the natural addition on $\Zzz_2$ and
$\Nat$ and the identical involution. The element  $u=(1,0)$ is like
in (\ref{Case1}).

Let us now consider the semigroup  $T=\Zzz_2\times\Nat$, with the
natural addition on $\Zzz_2$ and maximum as the operation on $\Nat$,
the involution is again set to identity. It is easy to see that $T$
is $*$--separable. The element $u=(1,0)$ is such that (\ref{Case2})
is satisfied. This example shows one more thing. Namely, the
condition
\begin{itemize}
\item[]{ \it $\dim\ker(u_\phi+I_{\h^\phi})<\infty$ for all $\phi\in\M S $ of compact support}
\end{itemize}
is \underbar{not} equivalent to  any of the conditions of Theorem \ref{main}. Indeed, the
characters on $T$ form a discrete, enumerable set. If the mapping $\phi\in\M S $ is
compactly supported then it is finitely supported and consequently the space $\h^\phi$ is
finite dimensional. Hence, $u=(1,0)$ satisfies the above condition, but does not satisfy
(i).

Note that in Proposition \ref{krein} restricting to compactly
supported moment functions is possible because the function $\phi$
constructed in the proof (ii)$\To$(i) is supported by only one
character.
\end{exm}

\section{Functions with a finite number of negative squares}

The condition $\dim\ker( u_\phi+I_{\h^\phi})<\infty$ can be written also in the language of
negative squares. Precisely speaking, by {\it the number of negative squares of a mapping }
$\psi:S\to \Comp$ satisfying
\begin{equation}\label{symmetric}
\psi(s)=\cl{\psi(s^*)}, \qquad s\in S,
 \end{equation}
 we understand the maximum, taken over all numbers $N\in\Nat$ and
all sequences $s_0\dts s_N$, of the number of negative eigenvalues
of the symmetric matrix $\left(\psi(s_i+s_j^*)\right)_{i,j=0}^N$.
Note that if $[u]=[u^*]$ and $\phi\in\Pp S$ then the mapping
$\psi:=\phi(\cdot+u)$ satisfies (\ref{symmetric}). Indeed, take any
character $\sigma$. Then
$\sigma(s+u)=\sigma(s)\sigma(u)=\sigma(s)\sigma(u^*)=\sigma(s+u^*)$.
By the result of \cite{biss-sep} we get $\phi(s+u)=\phi(s+u^*)$.
Combining this with $\phi(t)=\cl{\phi(t^*)}$ ($t\in S$)
\cite[4.1.6]{BChR} proofs the claim.

\begin{prop}\label{(v)}
Let $S$ be a $*$--semigroup with zero and let $u\in S$ be such that $[2u]=0$ and $[u]=[u^*]$
and let $\phi\in\Pp S$. Then the number of negative squares of the mapping $\phi(\cdot+u)$
equals $\dim\ker(u_\phi+I)$. Consequently, conditions {\rm (i)--(iv)} of Theorem \ref{main}
are equivalent to each of the following:
\begin{itemize}
\item[(v)] the mapping $\phi(\cdot+u)$ has a finite number of negative squares for every $\phi\in\Pp S$;
\item[(vi)] the mapping $\phi(\cdot+u)$ has a finite number of negative squares for every $\phi\in\M S $.
\end{itemize}
\end{prop}

\begin{proof} (cf. \cite{berg,langerio,sasvari} for similar arguments) First let
us assume that $\dim\ker(u_\phi+I)=m\in\Nat$. Consider the indefinite inner product space
$(\h^\phi, \seq{u_\phi\,\cdot\,,\,\cdot\,})$. Since  $\Dd^\phi$ is dense in $\h^\phi$ we can
find elements $s_1\dts s_k\in S$ and vectors $\alpha^i=(\alpha^i_1\dts
\alpha^i_k)\in\Comp^k$ ($i=1\dts m$) such that the elements
 \begin{equation}\label{fi}
 f^i:=\sum_{j=1}^k\alpha^i_j K^\phi_{s_j}\qquad (i=1\dts m)
 \end{equation}
span an $m$-dimensional negative subspace (\cite[Theorem IX.1.4]{bognar}). Let
 \begin{equation}\label{A}
A:=\left(\phi(s_j+s_{j'}^*+u)\right)_{j,j'=1}^k\in\Comp^{k\times k}.
 \end{equation}
  Note that for $i,l=1\dts m$ we have
\begin{equation}\label{fitransf}
 \seq{A \alpha^l,\alpha^i}=
 \sum_{j,j'=1}^k\alpha^{l}_j\cl{\alpha^{i}_{j'}}K^\phi(s_j+u,s_{j'})
 =\seq{u_\phi f^l,f^i}.
 \end{equation}
Hence, the subspace $\lin\set{\alpha^1\dts\alpha^k}$ is a negative subspace of the
indefinite inner product space $(\Comp^m,\seq{A\cdot,\cdot})$. Since $f^1\dts f^m$  are
linearly independent, the vectors $\alpha^1\dts\alpha^m$ are linearly independent as well.
Therefore, the matrix $A$ has at least $m$ negative eigenvalues.

Now let us assume that for some choice of $s_1\dts s_k\in S$ the
matrix $A$ defined as in (\ref{A}) has $m$ negative eigenvalues.
Then there exists linearly independent vectors
$\alpha^i=(\alpha^i_1\dts \alpha^i_k)\in\Comp^k$ ($i=1\dts m$) such
that
 \begin{equation}
 \seq{A \alpha^l,\alpha^i}=\delta_{il}\lambda_i\qquad
  i,l=1\dts m,
  \end{equation}
with some $\lambda_1\dts\lambda_m<0$. We define $f_1\dts f_m$ as in
(\ref{fi}) (with the new meaning of $s_1\dts s_m$). Using the
calculation in (\ref{fitransf}) we get that the space
$\lin\set{f_1\dts f_m}$ is a negative subspace of $(\h^\phi,
\seq{u_\phi\cdot,\cdot})$. We show now that $f_1\dts f_m$ are
linearly independent. If $\sum_{i=1}^m \beta_i f_i=0$ for some
$\beta_1\dts\beta_m\in\Comp$ then, by (\ref{fitransf}),
 $$
\seq{A \sum_{i=1}^m\beta_i\alpha^i,\sum_{i=1}^m\beta_i\alpha^i}=\seq{u_\phi \sum_{i=1}^m
\beta_i f_i, \sum_{i=1}^m \beta_i f_i}=0.
 $$
But $A$ is strictly negative on $\lin\set{\alpha^1\dts\alpha^m}$ and $\alpha_1\dts\alpha_m$
are linearly independent. Hence, $\beta_1=\cdots=\beta_m=0$.
\end{proof}

\section{More examples}
The reader can easily check that Theorem ~\ref{main} can be applied to Examples ~\ref{exm1}
and~\ref{exm2}. The next example concerns the estimation of the dimension of the eigenspace
in Theorem ~\ref{main}. We will show that this dimension can be any number between $0$ and
$M/2$, where $M$ is defined as in the proof of Theorem 2.

\begin{exm}
Let $S=\Zzz_2^m$ with identical involution and let $u=(1,0,0\dts0)$.
We have $M=2^m$. The dual semigroup $S^*$ can be identified with
$\set{-1,1}^m$. There are $2^{m-1}$ characters $\sigma$ on $S$
satisfying $\sigma(u)=-1$ and $2^{m-1}$ characters $\sigma$
satisfying $\sigma(u)=1$. We denote those characters by
$\sigma_1\dts\sigma_{2^{m-1}}$ and $\rho_1\dts\rho_{2^{m-1}}$,
respectively. For fixed $k,l\in\set{0\dts 2^{m-1}}$ we
put\,\footnote{\ We use the convention $\sum_{i=1}^0 a_i:=0$ }
$\mu:=\sum_{i=1}^k \delta_{\sigma_i}+\sum_{j=1}^{l}\delta_{\rho_j}$
and $\phi:=\Ll(\mu)$. Since the support of the measure is consists
of $k+l$ points, the space $\h^\phi$ is $k+l$ dimensional. To see
this one can use the interpretation of $\h^\phi$ as $\mathcal
P^\mu$, as in Example \ref{exm1}. Now let us observe that
\begin{equation}\label{sigmaker}
    \sigma_1\dts\sigma_k\in\ker( u_\phi+I),\quad
    \rho_1\dts\rho_{l}\in\ker( u_\phi-I),
\end{equation}
 by the same
argument as in the proof of Theorem \ref{main} (iii)$\To$(ii). Since the characters are
always linearly independent  we get that $\dim\ker( u_\phi+I)\geq k$ and $\dim\ker(
u_\phi-I)\geq l$. But the eigenspaces corresponding to $-1$ and $1$ are orthogonal, thus
$\dim\ker( u_\phi+I)=k$ and $\dim\ker( u_\phi-I)=l$.

Let us put $e=(0,1,0\dts 0)$ and take two numbers $l_1\in\set{0\dts 2^{m-1}}$ and
$l_2\in\set{0\dts 2^{m-2}}$. Using the same technique we can construct a mapping $\phi\in\Pp
S$ such that $\dim\ker(u_\phi+I)=l_1$ and $\dim\ker(e_\phi+I)=l_2$.
\end{exm}

In the following example there are three elements satisfying
$2u=0$ and $u=u^*$, with three different upper bounds for the
dimensions of the kernel.

\begin{exm}
 Let us consider the semigroup $S=\mathrm{U}(\Zzz_2^2,\Zzz_2,\pi)$ where
$\pi(x,y)=x$ for $x,y\in\Zzz_2$. The involution on $S$ is identity.  Note that $(1,0)+s\neq
s$ for $s\in S$, but $(0,1)+s\neq s$ only for $s\in\Zzz_2^2$. Hence, the upper bounds for
the dimensions of the kernels are three and two, respectively. The dimension of the kernel
for $(0,0)$ is obviously zero.
\end{exm}

\begin{rem}
Let us take two $*$-separative semigroups $S$ and $T$, both having  neutral elements ($0_S$
and $0_T$ respectively) and a $*$-homomorphism $h:S\to T$ satisfying $h(0_S)=0_T$. Take an
element $u\in T$ satisfying $2u=0_T$ and $u=u^*$. The $*$-semigroup $\mathrm{U}(S,T,h)$ has
a zero, namely $0_S$. However, the element $u$, understood as an element of
$\mathrm{U}(S,T,h)$, does not satisfy the condition $2u=0_{\mathrm{U}(S,T,h)}$.
Nevertheless, we have $3u=u$ and $u^*=u$, which by \cite{szaf-bdd} guaranties boundedness
 (and hence selfadjointness) of $u_\phi$ for any $\phi\in\Pp{\mathrm{U}(S,T,h)}$. The
indefinite inner product space $(\h^\phi,\seq{u_\phi\cdot,-})$ is then degenerate, i.e.
$u_\phi$ has a non-trivial kernel.
\end{rem}

 We could investigate, instead of positive definite mappings on
$S$, the set of positive definite forms on $S$. Namely, we say that $\phi:S\times
\mathcal{E}\times \mathcal{E}\to\Comp$ is a { form over ($S,\mathcal{E}$)} if for every
$s\in S$ the mapping $\phi(s,\,\cdot\,,-)$ is a hermitian bilinear form on the linear space
$\mathcal{E}$. We say that a form is positive definite if for every  finite sequences
$(s_k)_k\sbs S$, $(f_k)_k\sbs \mathcal{E}$ we have $\sum_{i,j} \phi(s_j^*+s_i;f_i,f_j))\geq
0$. For a positive definite form $\phi$ we can construct a Hilbert space $\h^\phi$ which
together with the functions $K_{s,f}^\phi$ ($s\in S$, $f\in\mathcal E$) constitute a RKHS.
Like in the case of $\mathcal{E}=\Comp$, cf. \cite{szaf-bdd} and also \cite{general}, we can
define the (closed) shift operator associated with an element $u\in S$ by
$u_\phi(K_{s,f}^\phi)=K_{s+u,f}^\phi$. %
%
%
The following example shows, that in this setting the equivalence
in Theorem \ref{main} no longer holds.
\begin{exm}

Let $S=\Zzz_2$ (with the identical involution) and let
 $\mathcal{E}$ be an \underbar{infinite} dimensional {Hilbert} space.
Consider the following form
    $$
    \phi(x,f,g)=\left\{\begin{array}{rcl}
                              \seq{f,g}_\mathcal{E} & : & x=0\\
                              \seq{-f,g}_\mathcal{E} & : & x=1
                              \end{array}\right. .
    $$
Note that
    $$
    \sum_{x,y=0,1}
    \phi(x+y,f_x,f_y)=\seq{f_0,f_0}+\seq{f_1,f_1}-2\Re\seq{f_1,f_0}=\norm{f_1-f_0}^2
    $$
which is greater or equal to zero for any choice of $f_0,f_1\in\h$.

The space $\h^\phi$ can be realized as $\h^\phi=\mathcal{E}$ so as $K_{0,f}^\phi=f$ and
$K_{1,g}=-g$, $f,g\in \mathcal{E}$.

Take $u=1$. The condition (i) of Theorem \ref{main} is satisfied
because the semigroup is of finite cardinality. On the other hand
 $u_\phi K_{0,f}^\phi = K_{1,f}^\phi=-f$ for $f\in\h$. Hence,
$\dim\ker(u_\phi+I)=\dim\h=\infty$.

\end{exm}

 \section{Final remarks}
 Our work is connected in a
  way with many other papers and books. Let us mention some of
them.
\begin{itemize}
 \item The transformation $\phi\mapsto\phi(\cdot+u)$ has been
 investigated by Bisgaard. He showed in \cite{biss-fact} that it
 is always a sum of four positive definite mappings. \item
 Functions with finite number of negative spaces on topological
 groups has been considered in the book of Sasv\'ari
 \cite{sasvari}. \item In \cite{berg}  sequences on
 $\Nat$ with a finite number of negative squares are considered.
 \item In \cite{BChR} the authors consider negative definite
 sequences, which is a subclass of mappings with a finite number
 of negative squares. \item Finally, in \cite{biss-cor}
 definitizing ideals are investigated.
 \end{itemize}

\bibliographystyle{amsalpha}

\begin{thebibliography}{99}




\bibitem{berg}{C. Berg, J.P.R. Christensen, P.H. Maserick, Sequences with finitely many negative
squares, {\it J. Funct. Anal.} {\bf79} (1988), 260--287.}




\bibitem{BChR}{Ch. Berg, J.P.R. Christensen, P. Ressel, {\it Harmonic Analysis on Semigroups,}
Springer-Verlag, New-York Berlin Heidelberg Tokyo {\bf1984}.}

\bibitem{bognar}{ J. Bogn\'ar, {\it Indefinite Inner Product
Spaces}, Springer-Verlag {\bf1974}.}



\bibitem{biss-scmat}{T.M. Bisgaard, On the Relation between the
Scalar Moment Problem and the Matrix Moment Problem on $*$-semigroups, \it Semigroup Forum
\rm {\bf68} (2005), 25-46. }

\bibitem{biss-two}{T.M. Bisgaard, Two sided complex moment problem,
{\it Ark. Mat.}  {\bf27}  (1989),   23--28.}

\bibitem{biss-sep}{T.M. Bisgaard, Separation by characters or positive definite functions, {\it Semigroup Forum} {\bf53} (1996), 317-320.}



\bibitem{biss-ext}{T.M. Bisgaard, Extensions of Hamburger's Theorem, {\it Semigroup Forum} {\bf57} (1998), 397-429.}



\bibitem{biss-CC}{T.M. Bisgaard, Semiperfect countable $\Comp$-separative $C$-finite semigroups, {\it Collect. Math.} {\bf52} (2001), 55-73}


\bibitem{biss-fact}{T.M. Bisgaard, Factoring of Positive Definite Functions
on Semigroups, {\it Semigroup Forum}  {\bf64} (2002), 243-264.}

\bibitem{biss-cor}{T.M. Bisgaard, H. Cornean,  Nonexistence in General of a Definitizing Ideal of the
Desired Codimension, {\it Positivity } {\bf7} (2003), 297–302.}


\bibitem{lacunae}{D. Cicho\'n, J. Stochel, F.H. Szafraniec, Extending positive
definiteness, preprint}


\bibitem{clipres}{A.H. Clifford, G.B. Preston, \it The Algebraic Theory of semigroups \rm
2 vols., AMS, Providence, {\bf1961, 1967}.}







\bibitem{langerio}{I.S. Iohvidov, M.G. Krein, H. Langer, {\it Introduction to the Spectral
 Theory of operators in Spaces with an Indefinite Metric},
 Akademie-Verlag, Berlin {\bf1982}.}







\bibitem{sasvari}{Z. Sasv\'ari, {\it Positive Definite and definitizable functions},
Akademie Verlag, Berlin {\bf1994}.}





\bibitem{szaf-bdd}{F.H. Szafraniec, Boundness of the shift operator
related to positive definite forms: an application to moment problems, {\it Ark. Mat.} {\bf
19} (1981), 251-259.}

\bibitem{test} F.H. Szafraniec,  Interpolation and domination by positive
definite kernels, in {\em Complex Analysis - Fifth Romanian-Finish Seminar}, Part 2,
Proceedings, Bucarest (Romania), 1981, eds. C. Andrean Cazacu, N. Boboc, M. Jurchescu and I.
Suciu, Lecture Notes in Math., {\bf 1014},  291-295, Springer, Berlin-Heidelberg, 1893.
   \bibitem{general} F.H. Szafraniec,  The Sz.-Nagy "th\'eor\`eme principal" extended. Application to subnormality,
   {\it Acta Sci. Math.  {\rm(}Szeged{\rm)} }, {\bf57}(1993), 249-262.

 \bibitem{test2} F.H. Szafraniec,   The reproducing kernel Hilbert space and its multiplication operators,
 {\em Oper. Theory Adv. Appl.}, {\bf114}(2000), 253-263.
%

\bibitem{nagy} B. Sz.-Nagy,     {\it Extensions   of
linear transformations in Hilbert space which extend beyond this space}, Appendix to F.
Riesz, B. Sz.--Nagy, {\it Functional Analysis}, Ungar, New York, {\bf1960}.



\bibitem{thill}{M. Thill, Exponentially bounded indefinite functions, Math. Ann. {\bf285} (1989),
297-307.}

\end{thebibliography}

\end{document}